\documentclass[12pt]{amsart}
\usepackage{amssymb,amsmath}
\makeatletter
\def\@cite#1#2{{\m@th\upshape\bfseries%
[{#1\if@tempswa{\m@th\upshape\mdseries, #2}\fi}]}}
\makeatother
\theoremstyle{plain}
\newtheorem{thm}{Theorem}[section]
\newtheorem{lem}[thm]{Lemma}
\newtheorem{cor}[thm]{Corollary}
\newtheorem{prop}[thm]{Proposition}
\theoremstyle{definition}
\newtheorem{rem}[thm]{Remark}
\newtheorem{defn}[thm]{Definition}
\newtheorem{note}[thm]{Note}
\newtheorem{eg}[thm]{Example}

\newcommand{\Prf}{\noindent\textbf{Proof.\ }}
\newcommand{\bx}{\strut\hfill$\blacksquare$\medbreak}

\newcommand{\ca}{\mathrm{C}^*}

\newcommand{\ol}{\overline}

\newcommand{\wot}{\textsc{wot}}

\newcommand{\bbC}{{\mathbb{C}}}

\newcommand{\bbF}{{\mathbb{F}}}
\newcommand{\bbN}{{\mathbb{N}}}

\newcommand{\bbZ}{{\mathbb{Z}}}
 \newcommand{\B}{{\mathcal{B}}}

 \newcommand{\E}{{\mathcal{E}}}
 \newcommand{\F}{{\mathcal{F}}}
 
\renewcommand{\H}{{\mathcal{H}}}

 \newcommand{\K}{{\mathcal{K}}}  

 \newcommand{\M}{{\mathcal{M}}}
 
\renewcommand{\O}{{\mathcal{O}}}

\renewcommand{\phi}{\varphi}
\newcommand{\upchi}{{\raise.35ex\hbox{$\chi$}}}
\newcommand{\fA}{{\mathfrak{A}}}
\newcommand{\fB}{{\mathfrak{B}}}

\newcommand{\qand}{\quad\text{and}\quad}

\newcommand{\qfor}{\quad\text{for}\quad}
\newcommand{\qwith}{\quad\text{with}\quad}

\newcommand{\spn}{\operatorname{span}}

\newcommand{\rowt}{(T_1, \ldots, T_N)}
\newcommand{\rows}{(S_1, \ldots, S_N)}
\newcommand{\rowl}{(L_1, \ldots, L_N)}

\newcommand{\fnplus}{{\bbF}_N^+}
\newcommand{\fock}{\ell^2 ({\bbF}_N^+)}
\newcommand{\alphiw}{\lambda_{i,w}}
\newcommand{\lamiw}{\lambda_{i,w}}
\newcommand{\cstarperk}{\mathrm{C}^*_N(\mathrm{per}\,\,k)}

\begin{document}

\title[Periodic Weighted Shifts on Fock Space]%
{Inductive Limit Algebras from Periodic Weighted Shifts on Fock space}
%
\author[D.W.Kribs]{David~W.~Kribs${}^1$} 
\address{Department of Mathematics, University of Iowa, Iowa City, IA 
52242}
\email{dkribs@math.uiowa.edu}
\thanks{${}^1$partially supported by a Canadian NSERC
Post-doctoral Fellowship.}  

\date{}
\begin{abstract}
Non-commutative multivariable versions of weighted shift operators arise
naturally as `weighted' left creation
operators acting on the Fock space Hilbert space. We identify a  natural 
notion
of periodicity for these $N$-tuples, and then  find a family of 
inductive limit 
algebras determined by the periodic 
weighted shifts which can be regarded as non-commutative 
multivariable generalizations of 
the Bunce-Deddens $\ca$-algebras. We establish this by proving that the
$\ca$-algebras generated by shifts of a given period are isomorphic to
full matrix algebras over Cuntz-Toeplitz algebras. This leads to
an isomorphism theorem which parallels the Bunce-Deddens and UHF  
classification scheme. 
\end{abstract}
\maketitle
\tableofcontents

The primary goal of this paper is to initiate the study of {\it
non-commutative multivariable
weighted shifts}. 
Almost three decades ago, Bunce and Deddens \cite{BD1,BD2} introduced a family of
inductive limit $\ca$-algebras generated by periodic unilateral weighted
shift operators. On the other hand, we now know that   
non-commutative multivariable versions of
unilateral shifts arise in theoretical physics
and free probability
theory as the so-called left creation operators acting on
the full Fock space Hilbert space. There is now an extensive body of 
research for these operators and the
algebras they generate (see \cite{AP,DKP,DP1,DP2,Kribs,Pop_fact,Pop_beur}
for example).
 
In this paper, we introduce a family of $\ca$-algebras which can
be regarded as non-commutative multivariable generalizations of the 
Bunce-Deddens
algebras. In
accomplishing this, based on the creation operators, we introduce the concept of a 
non-commutative multivariable weighted shift and discover a satisfying
notion of periodicity based on the structure of Fock space. 
We characterize these algebras in terms of inductive limits of full matrix
algebras over the
Cuntz-Toeplitz and Cuntz algebras. This leads to a classification theorem
which 
parallels the classification of UHF algebras by Glimm \cite{Glimm1}, 
and the          
Bunce-Deddens algebras classification \cite{BD1,BD2}, by supernatural
numbers. 

In the opening section we  recall the formulation of Fock space and the creation operators.
We also quickly review the basics of the Cuntz and Cuntz-Toeplitz 
algebras. In the second
section we introduce non-commutative weighted shifts and investigate
their basic structure.  The third section describes a pictorial method 
for thinking of these shifts, by using
the Fock space `tree' structure. This leads to a natural notion of periodicity, and then we define
the $\ca$-algebras we study in the rest of the paper. 
The final two sections consist of an in-depth analysis of these algebras. 
Most importantly, we prove they are 
isomorphic to inductive limits of full matrix algebras of 
distinguished sizes over 
Cuntz and Cuntz-Toeplitz 
algebras. Using this characterization we establish a classification 
theorem based on $K$-theory for the Cuntz 
algebras.

\section{Introduction}\label{S:intro}

We begin by recalling the formulation of the full Fock space Hilbert space 
and its associated creation operators. 
For $N \geq 2$, let $\fnplus$ be the unital free semigroup on $N$ {\it non-commuting} letters
$\{ 1, 2, \ldots, N \}$. We denote the unit in $\fnplus$ by $e$. One way 
to realize  
$N$-variable Fock space is as $\H_N = \fock$. From this point of view, the
vectors 
$\{ \xi_w : w\in\fnplus \}$ form an orthonormal basis for $\H_N$ which can
be thought of as a
generalized Fourier basis. The left creation operators (also known as the 
Cuntz-Toeplitz
isometries we will see below) $L = \rowl$ are defined on   $\H_N$ by their 
actions on basis vectors,  
\[
L_i \xi_w = \xi_{iw}, \qfor 1\leq i \leq N \qand w\in \fnplus.
\]
The $L_i$ are isometries with pairwise orthogonal ranges for which the sum
of the range projections
satisfies
$
\sum_{i=1}^N L_iL_i^* = I - P_e, 
$
where $P_e = \xi_e \xi_e^*$ is the rank one projection onto the span of the vacuum vector
$\xi_e$. We will discuss a helpful pictorial method for thinking of the actions of these operators
in Section 3.  

Most importantly for our purposes, this $N$-tuple forms the
non-commutative multivariable
version of a unilateral shift. This claim is well supported by a number of facts. 
For instance, the unilateral shift {\it is} obtained for $N =1$, and 
otherwise each of the $L_i$
is unitarily equivalent to a shift of infinite multiplicity. Further, the study of the $L_i$ in
operator theory and operator algebras was at least partly initiated by the dilation theorem of 
Frazho \cite{Fra1}, Bunce \cite{Bun} and Popescu \cite{Pop_diln}, which provided the 
non-commutative multivariable version of Sz.-Nagy's classical minimal
isometric dilation of a
contraction \cite{SF}. Namely, every row contraction of operators on
Hilbert space has
a minimal joint dilation to isometries, acting on a larger space, with pairwise orthogonal
ranges. The classical Wold decomposition shows that every isometry breaks up into an
orthogonal direct sum of a unitary together with copies of the shift. Analogously, Popescu's
version \cite{Pop_diln} shows that every $N$-tuple of
isometries with pairwise orthogonal ranges decomposes into an orthogonal direct sum of
isometries which form a representation of the Cuntz $\ca$-algebra $\O_N$ (see below), together
with copies of $L = \rowl$. 

In addition, the $\wot$-closed non-selfadjoint algebras generated by the
$L_i$ have been
shown by Davidson, Pitts, Arias, Popescu and others to be the appropriate 
{\it non-commutative analytic Toeplitz algebras} (see
\cite{AP,DP1,DP2,Kribs,Pop_beur}).
We also mention that the $\wot$-closed non-selfadjoint algebras generated
by the weighted shifts discussed here have been investigated in
\cite{Kribs_nonsa}, where a number of results from the single variable 
setting have been generalized, at the same time exposing new 
non-commutative phenomena.   
Finally, we note that compressing the creation operators to symmetric Fock space yields the 
{\it commutative} multivariable  shift. The $\ca$-algebras generated by 
weighted versions of which were
studied in \cite{CM} for instance. 
     
The $\ca$-algebras determined by the isometries $L = \rowl$ have also been
studied extensively. 
The $\ca$-algebra generated by $L_1, \ldots ,L_N$ is called the {\it
Cuntz-Toeplitz algebra} and is
denoted $\E_N$. The ideal generated by the rank one projection $ I -
\sum_{i=1}^N L_iL_i^*$ in
$\E_N$ yields a copy of the compact operators. When this ideal is factored out, the
$\ca$-algebra obtained is the {\it Cuntz algebra} $\O_N$. It is the universal 
$\ca$-algebra generated by the relations 
\[
s_i^* s_j = \delta_{ij} {\bf 1} \qfor 1\leq i,j \leq N \qand \sum_{i=1}^N
s_is_i^* = {\bf 1}.
\]
Up to isomorphism, $\O_N$ is the $\ca$-algebra generated by {\it any} $N$ isometries 
$S = \rows$ which satisfy these relations, since it is simple.

The $K$-theory for a $\ca$-algebra consists of a series of invariants which hold information on
equivalence classes of  projections in the matrix algebras over the algebra. 
The $K$-theory for $\O_N$ was worked out by Cuntz \cite{Cun_K}. In 
particular, its $K_0$ group is the finite abelian group  $K_0 (\O_N) = 
\bbZ / (N-1)\bbZ$.
In connection with classification results for inductive limits of Cuntz
algebras we mention work of Rordam \cite{Rordam}. We also note that our
isomorphism theorem has overlap with work of
Evans \cite{Evans}.

\section{Non-commutative Weighted Shifts}\label{S:wtd_shifts}

From the discussion in the previous section, we are led to the following 
definition for non-commutative multivariable weighted shifts. We shall 
drop the multivariable reference for succinctness. 
We mention that the idea for considering these weighted shifts came during the
author's preparation of \cite{Kribs_curv}, where a related class of $N$-tuples
was used in the analysis there.

\begin{defn}                                                
We say that an $N$-tuple of operators $S = \rows$ acting on a Hilbert space $\H$ forms a 
{\it non-commutative weighted shift} if there is a  unitary $U
: \H_N \rightarrow \H$, 
operators $T = \rowt$ on $\H_N$, and scalars $\{ \lambda_{i,w} : 1\leq i 
\leq N \qand w\in\fnplus\}$ such that $S_i = U T_i U^*$ for $1 \leq i 
\leq N$ and 
\[
T_i \xi_w = \lambda_{i,w} \xi_{iw} \qfor 1\leq i\leq N \qand w\in\fnplus.
\]
\end{defn}

\begin{note}
For the sake of brevity, we  assume that the weighted
shifts $T = \rowt$ we consider actually act on  Fock space
$\H_N = \fock$. 
Further, the proposition below will allow us to make the following 
simplifying assumption on weights throughout the paper:
\[
\mbox{{\it Assumption:       }} \lambda_{i,w} \geq 0 \qfor 1 \leq i \leq N
\qand w \in \bbF_N^+. 
\]
Indeed, every shift is jointly unitarily equivalent to a shift with 
nonnegative weights. 
\end{note}

\begin{prop}
Suppose $T=\rowt$ is a weighted shift with weights $\{\lambda_{i,w}\}$. 
Then 
there is a unitary $U\in\B(\H_N)$, which is diagonal with respect to the 
standard basis for $\H_N$, such that the weighted shift 
\[
(UT_1U^*,\ldots, UT_NU^*)
\]
has weights $\{|\lambda_{i,w}|\}$. 
\end{prop}

\Prf
We build the unitary by inductively choosing scalars $\mu_w$ and defining 
$U\xi_w = \mu_w\xi_w$. Put $\mu_e=1$. Let $k\geq 1$ and assume the scalars 
$\{\mu_w:|w|=k-1\}$ corresponding to words of length $k-1$ in $\bbF_N^+$ 
(as the empty word, the unit $e$ is taken to have length zero) have been 
chosen. The scalars $\{\mu_w:|w|=k\}$ are obtained in the following 
manner. For $iw\in\fnplus$ with $|w|=k-1$ and $1\leq i \leq N$, choose 
$\mu_{iw}\in\bbC$ of modulus one such that 
\[
(\ol{\mu_w}\,\lambda_{i,w})\,\mu_{iw} \geq 0.
\]
Now if $1\leq i \leq N$ and $w\in\fnplus$ are arbitrary, we have
\begin{eqnarray*}
\big( UT_iU^* \big) \, \xi_w = \ol{\mu_w} \, UT_i  \xi_w &=& 
\ol{\mu_w}\lambda_{i,w} U \xi_{iw} \\
&=& \big( \ol{\mu_w} \lambda_{i,w} \mu_{iw} \big) \, \xi_{iw}.
\end{eqnarray*}
This yields the desired conclusion.
\bx

We next present a direct generalization of the factorization of weighted 
shift operators into  products of the unilateral shift and  
diagonal weight 
operators.

\begin{prop}\label{norms}
Let $T = \rowt$ be a weighted shift. Then  each $T_i$  factors as $T_i = 
L_i W_i$, 
where $W_i$ is a positive operator
which is diagonal with respect to the standard basis for $\H_N$. It
follows that the norms of the
$T_i$ and the row matrix $T$ are given by
\begin{itemize}
\item[$(i)$] $|| T_i || = \sup \{ \alphiw \, :\, w\in\fnplus \}$ for $1\leq i \leq N$
\item[$(ii)$] $|| T || = \sup_{1\leq i \leq N} || T_i || = \sup \{ \alphiw \, :\, w\in\fnplus \,\,\, {\rm
and}\,\,\, 1\leq i \leq N \}$.
\end{itemize}
\end{prop}

\Prf
For $1\leq i \leq N$, the operators $W_i$ are given by the equation 
\[
W_i \xi_w = (T_i\xi_w , \xi_{iw} ) \xi_w  = \alphiw \xi_w. 
\]
Since $W_i \geq 0$, we have $T_i^*T_i = W_i^*L_i^* L_i W_i = W_i^2$,  
which is diagonal.
Hence, $W_i = (T_i^* T_i)^{1 / 2}$ and $T_i = L_i W_i$. 
Further, this shows that 
\begin{eqnarray*}
|| T_i ||^2 = || T_i^* T_i || = || W_i^2 || &=& \sup \{ ||W_i^2 \xi_w || : 
w\in\fnplus \} \\ 
&=& \sup \{ \alphiw^2 : w\in\fnplus
\}. 
\end{eqnarray*}

On the other hand, the entries of the $N \times N$ matrix $T^*T$ consist of $T_i^*T_i$'s down
the diagonal and zero off the diagonal, since the ranges of the $T_i$ are pairwise orthogonal.
Hence, from the above computation
\[
||T|| = ||T^*T||^{1 / 2} = \sup_{1\leq i \leq N} || T_i ||  = \sup_{1\leq i \leq N} \{ \alphiw : w \in
\fnplus \},
\]
which completes the proof. 
\bx

We finish this section by  observing that the $\ca$-algebra generated by a
non-commutative weighted shift, which is bounded below in an appropriate sense, contains the
Cuntz-Toeplitz algebra.

\begin{defn}
Let $T = \rowt$ be a weighted shift. If each $W_i$ is bounded away from 
zero, in other
words if
\[
\inf \{ \alphiw : 1\leq i\leq N \,\,\, {\rm and} \,\,\, w \in \fnplus \} > 0,
\]
we say that $T$ is {\it bounded below}.
\end{defn}

\begin{cor}\label{bded_below}
The $\ca$-algebra $\ca (T_1,\ldots, T_N)$ generated by the operators $\{ 
T_1, \ldots, T_N\} $ from a
weighted shift
$T = (T_1, \ldots ,T_N)$ contains $L = \rowl$ 
when $T$ is  bounded below.
\end{cor}

\Prf
From the proof of the previous proposition, we see that $W_i$ is invertible precisely when 
$\inf \{ \alphiw : w \in \fnplus \} > 0$. Thus, $T$ being bounded below implies that each
$W_i$ is invertible. However, $W_i = (T_i^*T_i)^{1 / 2}$ belongs to $\ca (T_i)$, and hence to
$\ca (T_1, \ldots , T_N)$, thus so does $L_i = T_i W_i^{-1}$ for $1 \leq i \leq N$. 
\bx

\begin{note}
We mention that the $\ca$-algebras $\ca(T_1,\ldots, T_N)$ generated by the 
$T_i$ from a single
weighted shift $T
= \rowt$ are the focus of
analysis in \cite{CKM}. 
\end{note}

\section{Fock Space Trees And Periodicity}\label{s:trees}

In this section we aim to convey to the reader a helpful pictorial method 
for thinking of 
non-commutative weighted shifts. In doing so, we introduce what seems to 
be a natural notion
of periodicity for these operators. We also define the operator algebras which will be studied in
the rest of the paper. 

Recall that $N$-variable Fock space $\H_N = \ell^2 (\bbF_N^+)$ has 
the orthonormal basis $\{ \xi_w : w\in\fnplus\}$. This basis
yields a natural tree
structure for Fock space which is traced out by the creation operators, 
and more generally by 
weighted shifts.  

\begin{defn}
Let $T = \rowt$ be a  weighted shift. 
Let $\F_T$ be the set of vertices $\{ w : w\in\fnplus \}$, together with the `weighted' directed
edges which correspond to the directions 
\[
\{ \lambda_{i,w} := w \mapsto iw\,\,\, \big|  \,\,\, {\rm for} \,\,\, 1 \leq i \leq N \,\,\, {\rm and} \,\,\,
w \in\fnplus\}.
\]
We regard an edge $\lambda_{i,w}$ as lying to the left of another edge $\lambda_{j,w}$
precisely when $i<j$. We call $\F_T$ the {\it weighted Fock space tree} 
generated by $T$. 
\end{defn}

Pictorially, with $N=2$ as an example, a typical weighted Fock space tree 
is given by the following diagram:

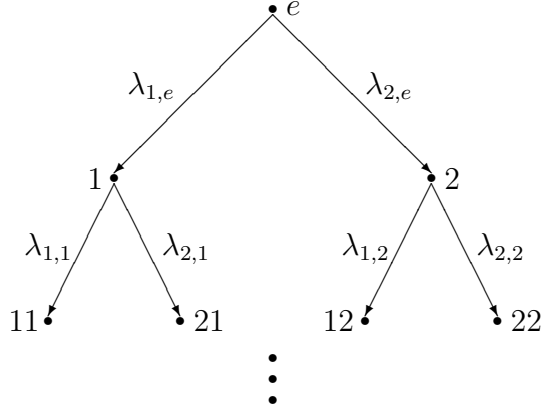
\begin{picture}(100,165)(-110,-26)
\put(50,132){\circle*{3}}
\put(-10,68){\circle*{3}}
\put(110,68){\circle*{3}}
\put(55,130){$e$}
\put(-20,64){1}
\put(115,64){2}
\put(50,130){\vector(-1,-1){60}}
\put(50,130){\vector(1,-1){60}}
\put(-5,100){$\lambda_{1,e}$}
\put(85,100){$\lambda_{2,e}$}
\put(-10,66){\vector(1,-2){25}}
\put(-10,66){\vector(-1,-2){25}}
\put(110,66){\vector(1,-2){25}}
\put(110,66){\vector(-1,-2){25}}
\put(-35,14){\circle*{3}}
\put(15,14){\circle*{3}}
\put(85,14){\circle*{3}}
\put(135,14){\circle*{3}}
\put(-50,10){11}
\put(20,10){21}
\put(69,10){12}
\put(140,10){22}
\put(50,0){\circle*{3}}
\put(50,-8){\circle*{3}}
\put(50,-16){\circle*{3}}
\put(-44,39){$\lambda_{1,1}$}
\put(76,39){$\lambda_{1,2}$}
\put(8,39){$\lambda_{2,1}$}
\put(127,39){$\lambda_{2,2}$}
\end{picture}

{\noindent}Notice that this structure is really determined by the operators $T = \rowt$.
Indeed, given a basis vector $\xi_w$, the directed edge $\lambda_{i,w}$ corresponds to the
action of $T_i$ on $\xi_w$, namely mapping it to $\lamiw \xi_{iw}$. 
Thus, more generally, we have the following picture for weighted edges leaving a typical vertex
in the tree:

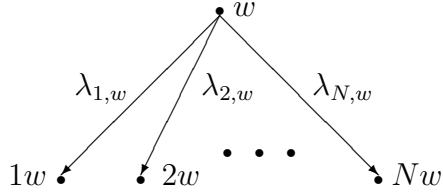
\begin{picture}(100,82)(-110,60)
\put(50,132){\circle*{3}}
\put(-10,68){\circle*{3}}
\put(110,68){\circle*{3}}
\put(55,130){$w$}
\put(-30,66){$1w$}
\put(20,68){\circle*{3}}
\put(28,66){$2w$}
\put(115,66){$Nw$}
\put(50,130){\vector(-1,-1){60}}
\put(50,130){\vector(1,-1){60}}
\put(-5,100){$\lambda_{1,w}$}
\put(43,100){$\lambda_{2,w}$}
\put(85,100){$\lambda_{N,w}$}
\put(50,130){\vector(-1,-2){30}}
\put(53,78){\circle*{3}}
\put(65,78){\circle*{3}}
\put(77,78){\circle*{3}}
\end{picture}

There are a number of conceptual benefits obtained by identifying these 
trees with weighted
shifts.  For instance, this point of view leads to the following notion of  
periodicity. 

\begin{defn}
Let $k\geq 1$ be a positive integer. We say that a 
weighted shift $T = \rowt$ is of {\it period $k$} if 
\[
T_i \xi_w = \lambda_{i,u} \xi_{iw} \qfor w\in\fnplus,
\]
where $w=uv$ is the {\it unique} decomposition of $w$ with $0\leq |u| < k$ 
and $|v| \equiv 0 ({\rm
mod}\,\,k)$. 
\end{defn}

\begin{note}
Observe that this says the scalars $\{ \lambda_{i,u} : 0\leq |u| < k \}$ 
{\it completely} determine
the shift. They can be thought of as a `remainder tree top'. 
For $N=1$ the standard notion of periodicity is recovered, 
since the tree collapses to a single infinite stalk.  For $N\geq 2$ It is 
most
satisfying to think of this notion of periodicity in terms of the tree 
structure: If 
$T = \rowt$ is period $k$, then the remainder tree top, that is the finite 
top of the tree determined by vertices $\{ w : |w| < k
\}$ and edges $\{ \alphiw : |w|<k \}$, is `repeated' throughout the entire 
weighted 
tree. 

In fact, this finite tree top is repeated with a certain exponential 
growth. For instance, at the
level of the tree corresponding to words $w\in\fnplus$ of length $nk$ for some positive integer
$n\geq 1$, the top of this finite tree is repeated $N^{nk}$ times, once for every word of
length $nk$. 

We mention that related tree top constructions play a key role in the 
paper \cite{KP}. 
\end{note}

Finally, we introduce the operator algebras which we are interested in 
studying. 

\begin{defn}
For positive integers $N \geq 2$ and $k\geq 1$,  let $\cstarperk$ be 
the $\ca$-algebra (contained in $\B
(\H_N)$)  generated by the $T_i$ from all  weighted shifts
$T = \rowt$ of period $k$.

It is clear from the picture given by the Fock space trees that if $n_1 | n_2$, then
$\mathrm{C}^*_N(\mathrm{per}\,\,n_1)$ is contained in 
$\mathrm{C}^*_N(\mathrm{per}\,\,n_2)$. (We prove this at the end
of this section.)
Thus, given  an increasing sequence of positive integers $\{ n_k \}_{k\geq 1}$ with $n_k |
n_{k+1}$ for $k \geq 1$, we may consider the  inductive limit algebra
\[
\fA (n_k) = \overline{\bigcup_{k\geq 1}  \mathrm{C}^*_N(\mathrm{per}\,\,n_k) }
\]
determined by this sequence. Let $q$ be the quotient map
of $\B(\H_N)$ onto the Calkin algebra. We are also interested in describing the inductive
limit algebras $q ( \fA(n_k) )$.
\end{defn}

\begin{note}
The reader may find it helpful to know that $\cstarperk$ is generated by
the $T_i$ from a single
weighted shift. This is proved in the next section, using the matrix
decompositions obtained there. 
\end{note}

\begin{prop}\label{contain}
Let $n_1, n_2$ be positive integers with $n_1 | n_2$. Then
$\mathrm{C}^*_N(\mathrm{per}\,\,n_1)$ is contained in 
$\mathrm{C}^*_N(\mathrm{per}\,\,n_2)$.
\end{prop}

\Prf
Let $T = \rowt$ be a period $n_1$ weighted shift. Let $w \in \bbF_N^+$ and write $w = u_2
v_2$ with $0 \leq |u_2| < n_2$ and $|v_2| \equiv 0 ({\rm mod}\,\,n_2)$. We must show that
$\lamiw = \lambda_{i,u_2}$ for $1 \leq i \leq N$. In other words, $T_i \xi_w = \lambda_{i,u_2}
\xi_{iw}$ for each $i$. To see this, write $u_2 = u_1 v_1$ where $0 \leq |u_1| < n_1$ and $|v_1|
\equiv 0 ({\rm mod}\,\,n_1)$. Since $T$ is of period $n_1$ we have
\[
T_i \xi_{u_2} =   \lambda_{i,u_2} \xi_{iu_2}  = \lambda_{i,u_1} \xi_{iu_2}, 
\]
so that $  \lambda_{i,u_2} =  \lambda_{i,u_1}$.

On the other hand, since $n_1 | n_2$ we have $|v_2| \equiv 0 ({\rm mod}\,\,n_1)$. This tells us
that  $w = u_1 (v_1 v_2)$ with $0 \leq |u_1| < n_1$ and $|v_1 v_2| \equiv 0 ({\rm mod}\,\,n_1)$.
Thus, $n_1$-periodicity once again gives us
\[
T_i \xi_{w} = \lambda_{i,w} \xi_{iw} =  \lambda_{i,u_1} \xi_{iw},
\]
Hence $\lamiw =  \lambda_{i,u_1} =  \lambda_{i,u_2}$, as required. 
\bx

\section{Main Theorem}\label{S:simple}

The $\ca$-algebra $\cstarperk$ generated by the $k$-periodic weighted
shifts can be
described in terms of a full matrix algebra with entries in a  
Cuntz-Toeplitz
algebra. From the discussion in Section 1, recall the Cuntz-Toeplitz 
algebra
$\E_{N^k}$ is the $\ca$-algebra generated by the creation operators $L =
(L_1, \ldots, L_{N^k})$ acting on $N^k$-variable Fock space $\H_{N^k}$. 

\begin{thm}\label{periodk}
For positive integers $N \geq 2$ and $k \geq 1$, let $d_{N,k}$ be the total number of words
in $\fnplus$ of length strictly less than $k$; that is,  
$
d_{N,k} = 1 + N + \ldots + N^{k-1}.
$
Then the algebra $\cstarperk$ of $k$-periodic weighted shifts is unitarily
equivalent to the algebra $\M_{d_{N,k}} 
(\E_{N^k})$ of $d_{N,k}\times d_{N,k}$ matrices with entries in 
$\E_{N^k}$. Further, this algebra is
generated by the $T_i$ from a single shift $T = \rowt$.  
\end{thm}                                         

\begin{rem}
At first glance the $N^k$ appearing in the theorem may seem somewhat
peculiar to the reader. We shall see that it arises from the
exponential nature of  periodicity here. We mention that the special case 
$N=2$ and $k=2$ of the theorem is expanded on in Example~\ref{thmeg}. 
\end{rem}

We shall prove the theorem in several stages. Throughout, $N \geq 2$ and 
$k\geq 
1$ will be
fixed positive integers.  The first step is
to decompose Fock space in a manner
which will lead to simple
matrix representations of the  periodic weighted shifts. 

\begin{lem}\label{spatial_decomp}
For $w \in \fnplus$ with $|w| < k$, the subspaces $\K_w$ of $N$-variable 
Fock space $\H_N$
given by
\[
\K_w = \spn \{ \xi_{wv} : |v| = km, \,\, m\geq 0 \},
\]
are pairwise orthogonal and 
\[
\H_N = \sum_{|w| < k} \oplus \, \K_w.
\]
Further, for $|w| < k$, the operators $U_w : \K_e \rightarrow \K_w$ 
defined by $U_w \xi_{v} =
\xi_{wv}$, for $|v| = km$ with $m \geq 0$, are unitary. Hence 
\[
U := \sum_{|w| <k}\! \oplus \,U_w \,\, : \,\, \K_e^{(d_{N,k})} 
\longrightarrow \H_N
\]
is a unitary operator. 
\end{lem}
     
\Prf
The subspaces $\K_w$ for $|w| < k$ clearly span $\H_N = \fock$ since any 
word $u \in \fnplus$
can be written, in fact uniquely, as $u = u_1 u_2$, where $|u_1| < k$ and $k$ divides $|u_2|$. To
see orthogonality, let $w_1$, $w_2$ be words with $|w_i|< k$, and consider typical basis vectors 
$\xi_{w_i v_i}$ for $\K_{w_i}$, where $|v_i | = km_i$ and $m_i \geq 0$. 
The only way the
inner product
$( \xi_{w_1v_1} ,  \xi_{w_2v_2} )$ can be non-zero, is if $w_1v_1 = w_2v_2$. But then we
would have 
\[
|w_1| + km_1 = |w_1| + |v_1| = |w_2| + |v_2| = |w_2| + km_2,
\]
so that $|w_1| = |w_2| < k$ and $m_1 = m_2$. This would imply that $w_1 = w_2$ and  
$v_1 = v_2$. It follows from this calculation that the subspaces $\{\K_w : 
|w|<k\}$ are pairwise
orthogonal. 

The operators $U_w$ as defined are unitary since they send one orthonormal basis to another.
Spatially, these unitaries can be thought of as the restrictions of
the isometries $L_w$, where $L_w := L_{i_1} \cdots L_{i_s}$ when $w$ is
the word $w = i_1
\cdots i_s$, to a distinguished subspace $\K_e$ of Fock
space. Alternatively, the action of the
adjoint $U_w^*$ on $\K_w$ is described by restricting
$L_w^*$ to $\K_w$. The last statement of the lemma is immediate from the 
spatial
decomposition of $\H_N$. 
\bx

We will distinguish between coordinate spaces of $\K_e^{(d_{N,k})}$ in the 
following manner:
For $w\in\fnplus$ with $|w| < k$, let 
\[
\{\, \xi_u^w \, : \, u\in\fnplus \qwith |u| = km \qfor m\geq 0 \}
\]
be the standard basis for the {\it $w$th coordinate space} of 
$\K_e^{(d_{N,k})}$, which is
given by $U^* \K_w = U_w^* \K_w$. Notice that for $|v|, |w| < k$ and $|u| 
=km$, the vectors
$\xi_u^w$ and
$\xi_u^v$ really correspond to the {\it same} vector $\xi_u$ in $\K_e$. 
Further, the action of
$U$ is
described by 
\begin{eqnarray}\label{udefn}
U \xi^w_u = \xi_{wu},
\end{eqnarray}
for $w,u\in\bbF_N^+$ with $|w| <k$ and $|u| = km$. 

The next step is to point out a relationship between particular Fock space trees. 

\begin{defn}
We define a natural bijective correspondence between; the $N^k$ words of
length $k$ in
$\fnplus$ on the one hand, and the $N^k$ letters which generate $\bbF_{N^k}^+$ on the other,
through the function
\[
\phi: \{w\in\fnplus:|w|=k\} \longrightarrow \{ w\in\bbF_{N^k}^+:|w|=1\}
\]
given by
\[
\phi(i_1 i_2\cdots i_k) = (i_1-1)N^{k-1} + \ldots + (i_{k-1}-1)N + i_k,
\]
for $1\leq i_j \leq N$ and $1 \leq j \leq k$. This correspondence is also characterized by
associating the words $\{ iw\in\fnplus : |w| = k-1\}$ with the `$i$th block' of $N^{k-1}$ letters
in the listing $\{ 1,2,\ldots,N^k\}$. Notice with this ordering that the operators 
$\{ L_{\phi(w)} : w\in\bbF_N^+, |w| =k\}$ are the $N^k$ creation operators associated with 
$N^k$-variable Fock space $\H_{N^k}$. 

The map $\phi$ extends in a natural way to a bijective identification
$\phi_m$ of
the set $\{
w\in\fnplus: |w| =km\}$ with $\{ w\in\bbF_{N^k}^+ : |w| =m\}$ for $m\geq
0$. Given $w_1,\ldots, w_m \in\bbF_N^+$ with $|w_i|=k$, the
extensions are given by 
\[
\phi_m (w_1 \cdots w_m) = \phi(w_1 ) \cdots \phi(w_m), 
\]
The units in $\bbF_N^+$ and $\bbF_{N^k}^+$ are
identified with each other. We will use the notation $\phi$ for the
extended map as well.
This ordering leads to the following spatial equivalence. 
\end{defn}

\begin{lem}\label{spatial_equiv}
The map from $\K_e = \spn \{ \xi_w : w\in\fnplus, |w| =km, m\geq 0 \}$ to 
$N^k$-variable Fock space $\H_{N^k} = \ell^2(\bbF_{N^k}^+)$ which sends a basis
vector $\xi_{w_1 \cdots w_m}\in\K_e$, where each $|w_i| = k$, to the basis 
vector
$\xi_{\phi(w_1) \cdots \phi(w_m)}\in\H_{N^k}$, is unitary.
\end{lem}

\Prf
This follows directly from the definitions of the space $\K_e$ and the map 
$\phi$. 
\bx

This lemma gives us a tight spatial equivalence between the orthogonal direct sums 
$\K_e^{(d_{N,k})}$ and $\H_{N^k}^{(d_{N,k})}$, which carries through for 
the weighted
shifts. We wish to preserve the correspondence discussed after Lemma~\ref{spatial_decomp}.
In particular, for $w\in\fnplus$ with $|w| < k$, we
let 
\[
\{ \, \xi_{\phi(u)}^w \, : \, u\in\fnplus \qwith |u| = km \qfor m\geq 0 \}
\]
be the standard basis for the $w$th coordinate space of  $\H_{N^k}^{(d_{N,k})}$. Once again,
with this identification, for 
$|v|, |w| < k$ and $|u| =km$, the vectors $\xi_{\phi(u)}^w$ and $\xi_{\phi(u)}^v$ 
correspond to the same vector $\xi_{\phi(u)}$ in $\H_{N^k}$. Finally, we let $V : 
\H_{N^k}^{(d_{N,k})} \rightarrow \K_e^{(d_{N,k})}$ be the unitary operator 
which encodes
this correspondence and the action of the map from the previous lemma. For $w, u \in\fnplus$
with $|w| < k$ and $|u| = km$,  the action of $V$ is given by 
\begin{eqnarray}\label{vdefn}
V\xi_{\phi(u)}^w = \xi_{u}^w.
\end{eqnarray}

With these Fock space decompositions in hand, we are now ready to focus on
the particular actions of
weighted shifts. 

\begin{lem}\label{shift}
For $w\in\fnplus$ with $|w| < k$, let $P_w$ be the orthogonal projection of
$\H_{N^k}^{(d_{N,k})}$ onto the $w$th coordinate space of $\H_{N^k}^{(d_{N,k})}$, so
that $I = \sum_{|w|<k}\oplus\,P_w$. Let $T = \rowt$ be a $k$-periodic 
weighted shift acting on $\H_N$. Then the operators 
${\rm Ad}_{UV} (T_i) = V^*U^*T_iUV$ act on  $\H_{N^k}^{(d_{N,k})}$ and 
have the
following block matrix decompositions:
\begin{itemize}
\item[$(i)$] For $|w| < k-1$ and $|v| <k$,
\[
P_v \big( {\rm Ad}_{UV} (T_i)  \big) P_w = 
\left\{ \begin{array}{ll}
\lamiw I_{\H_{N^k}} & \mbox{if $v=iw$} \\
0         & \mbox{if $v\neq iw$}
\end{array}
\right. 
\]
\item[$(ii)$] For $|w| = k-1$ and $|v| <k$,
\[
P_v \big( {\rm Ad}_{UV} (T_i)  \big) P_w = 
\left\{ \begin{array}{ll}
\lamiw L_{\phi(iw)} & \mbox{if $v=e$} \\
0         & \mbox{if $v\neq e$}
\end{array}
\right. 
\]
\end{itemize}
\end{lem}

\Prf
We first prove case $(i)$. From the preceding discussion, the vectors $\xi^w_{\phi(u)}$, where 
$u \in\fnplus$ with $|u| = km$, form an orthonormal basis for the range 
of $P_w$. Further, from equations (\ref{udefn}) and (\ref{vdefn}) we have 
\[
UV \xi^w_{\phi(u)} = U\xi_u^w =  \xi_{wu}.
\] 
Lastly, as $T = \rowt$ is $k$-periodic, we have $T_i \xi_{wu} = \lamiw \xi_{iwu}$. Since 
$|w| < k-1$, these facts lead us to the following computation:
\begin{eqnarray*}
P_v \big( {\rm Ad}_{UV} (T_i) \big) \xi^w_{\phi(u)} &=& P_v V^* U^* T_i \xi_{wu} \\
   &=& \lamiw P_v V^* U^* \xi_{(iw)u} \\
   &=& \lamiw P_v \xi_{\phi(u)}^{iw} = \lamiw \delta_{v,iw} \, \xi_{\phi(u)}^{iw}, 
\end{eqnarray*}
where $\delta_{v,iw}$ is equal to 1 if $v = iw$, and is 0 otherwise.  But recall that the vectors  
$ \xi_{\phi(u)}^{w}$ and $\xi_{\phi(u)}^{iw}$ both correspond to
the same vector $\xi_{\phi(u)}$ in $\H_{N^k}$. Thus, case $(i)$ is established. 

Now suppose that $|w| = k-1$. Again, $k$-periodicity gives us $T_i \xi_{wu} = \lamiw
\xi_{iwu}$. In this case $|iw| = k$, hence the definition of $U$ and 
$V$ from (\ref{udefn}) and (\ref{vdefn}) yields:
\begin{eqnarray*}
P_v \big( {\rm Ad}_{UV} (T_i) \big) \xi^w_{\phi(u)} &=& P_v V^* U^* T_i \xi_{wu} \\
   &=& \lamiw P_v V^* U^* \xi_{(iw)u} \\
   &=& \lamiw P_v V^* \xi^e_{iwu} \\
   &=& \lamiw P_v \xi_{\phi(iwu)}^{e} = \lamiw \delta_{v,e} \, \xi_{\phi(iwu)}^{e}. 
\end{eqnarray*}
However, from the definition of $\phi$, we have $\phi(iwu) = \phi(iw) \phi(u)$, since $|iw| = k$.
Therefore,
\[
P_v \big( {\rm Ad}_{UV} (T_i) \big) \xi^w_{\phi(u)} = \lamiw \delta_{v,e} L_{\phi(iw)}
\xi^e_{\phi(u)}.
\]
Once again, the vectors $\xi^w_{\phi(u)}$ and $\xi^e_{\phi(u)}$ both correspond to
$\xi_{\phi(u)}$ in $\H_{N^k}$. This establishes case $(ii)$, and completes the proof.
\bx

{\noindent}{\bf Proof of Theorem~\ref{periodk}}
We define an injective homomorphism of $\ca$-algebras $\pi: \cstarperk \rightarrow
\M_{d_{N,k}}(\E_{N^k})$ by $\pi(T_i) = {\rm Ad}_{UV} (T_i)$, for
every $k$-periodic  weighted shift $T= \rowt$. The map $\pi$ is clearly
an injective homomorphism since it is a unitary equivalence.  Further, it follows from case $(i)$
in
Lemma~\ref{shift}, that all the matrix units in $\M_{d_{N,k}}(\E_{N^k})$ can be obtained in
the image of $\pi$,  by
judicious choice of scalars $\lamiw$'s and appropriate matrix multiplication.
From case $(ii)$
in that lemma, we see that the $N^k$ creation operators which generate $\E_{N^k}$ can be
obtained in certain matrix entries. Since all the matrix units are present in the image, these
creation operators can be moved around to every entry. Therefore, it follows that $\pi$ is also
surjective, and hence defines a $*$-isomorphism. 

Lastly, it is not hard to see from the matrix decompositions of Lemma~4.6 
that the algebra $\cstarperk$ is generated by the $T_i$ from a single 
shift $T= \rowt$. For instance, from work in \cite{CKM} it follows that 
any shift will do for which the $N^k$ numeric $k$-tuples corresponding to 
the weights on each path of length $k$ in the associated tree are 
different. 
\bx

Before continuing, we discuss a special case of the theorem which may help 
to clarify some of the technical issues. 

\begin{eg}\label{thmeg}
Consider the case when $N=2$ and $k=2$. Then $\ca_2({\rm per}\,2)$ is the 
$\ca$-algebra generated by the $T_i$ from all 2-periodic shifts 
$T=(T_1,T_2)$. The theorem shows that this algebra is unitarily equivalent 
to the matrix algebra $\M_3(\E_4)$. Let us expand on this point. 

Such 2-tuples act on the Fock space $\H_2$, which has orthonormal basis 
$\{\xi_w:w\in\bbF^+_2\}$. As in the previous section, the remainder tree 
top which determines the weighted Fock space tree for a given 2-periodic 
shift $T=(T_1,T_2)$ is generated by six scalars $\{a,b,c,d,e,f\}$ as 
follows:
\[ 
T_1\xi_e = a\,\xi_1, \,\,\, 
T_1 \xi_1 = c\,\xi_{1^2}, \,\,\,
T_1 \xi_2 = e\,\xi_{12}.
\]
and
\[
T_2\xi_e = b\,\xi_2, \,\,\,
T_2\xi_1 = d\,\xi_{21}, \,\,\,
T_2 \xi_2 = f\,\xi_{2^2}.
\]
Thus, for example, the action of $T_1$ on basis vectors is given by 
\[
T_1 \xi_w = 
\left\{ \begin{array}{ll}
a\,\xi_{1w} & \mbox{if $|w|$ is even}\\
c\,\xi_{1w} & \mbox{if $w=1v$ and $|v| $ is even}\\
e\,\xi_{1w} & \mbox{if $w=2v$ and $|v|$ is even}.
\end{array}\right.
\]

In the proof of the theorem for this case, 2-variable Fock space $\H_2 = 
\ell^2(\bbF_2^+)$ decomposes into a direct sum $\H_2 = \K_e \oplus \K_1 
\oplus \K_2$ of $d_{2,2} = 1+2=3$ subspaces, each of which may be 
naturally identified with $(2^2=)4$-variable Fock space $\H_4 = 
\ell^2(\bbF^+_4)$. Take $\K_e$ for example. It is given by
\[
\K_e = \spn\Big\{ \xi_e, \{\xi_{1^2},\xi_{12},\xi_{21},\xi_{2^2}\}, 
\{ \xi_w : w\in\bbF_2^+, |w|=4\}, \ldots \Big\}. 
\]

The unitary equivalence produced by this spatial identification yields the 
following block matrix form for our given 2-periodic shift $T=(T_1,T_2)$, 
with respect to the decomposition $\H_2 = \K_e\oplus\K_1\oplus\K_2 \simeq 
\H_4\oplus \H_4 \oplus \H_4$, 
\[
T_1 \simeq \left[
\begin{matrix}
0 & cL_1 & eL_3 \\
aI & 0 & 0 \\
0&0&0
\end{matrix}\right]
\,\,\,\,\,\,
{\rm and} 
\,\,\,\,\,\,\,\,
T_2 \simeq \left[ 
\begin{matrix}
0 & dL_2 & fL_4 \\
0&0&0 \\ 
bI & 0 &0 
\end{matrix}\right],
\]
where $\{L_1,L_2,L_3,L_4\}$ are the standard creation operators on $\H_4$. 

Since we have complete freedom in $\ca_2({\rm per}\,2)$ on the choices of 
scalars $\{a,b,c,d,e,f\}$, it is now easy to see why it is unitarily 
equivalent to the matrix algebra $\M_3(\E_4)$. Further, it follows from 
these matrix decompositions that $\ca_2({\rm per}\, 2)$ is generated, for 
instance, by $\{T_1,T_2\}$ with $a=b=1$, $c=1/2$, $d=1/4$, $e=1/8$, 
$f=1/16$.  
\end{eg}

From Theorem 4.1 it follows that 
when we factor out the ideal of compact operators from $\cstarperk$, 
simple $\ca$-algebras are obtained. The key
point
being that the Cuntz-Toeplitz algebra is the extension of the compacts by the Cuntz algebra. 

\begin{cor}\label{qperk}
Let $q$ be the quotient map of $\B(\H_N)$ onto the Calkin algebra. Then for $N\geq 2$ and
$k\geq 1$, the algebra $q ( \cstarperk )$ is $*$-isomorphic to the matrix algebra
$\M_{d_{N,k}} (\O_{N^k})$. In particular, it is a simple $\ca$-algebra. 
\end{cor}

It follows that the inductive limit algebras  $q ( \fA (n_k) )$  defined in the previous section are
simple and have real rank zero. 

\begin{cor}
Let $N\geq 2$ and let $\{ n_k \}_{k\geq 1}$ be an increasing sequence of positive integers such
that $n_k$ divides $n_{k+1}$ for $k \geq 1$. Then the inductive limit algebra 
$q ( \fA (n_k) )$ is simple and has real rank zero.
\end{cor}

\Prf
Every ideal in $q( \fA (n_k))$ is the closed union of ideals of the
subalgebras
$q (  \mathrm{C}^*_N(\mathrm{per}\,\,n_k)  )$.
Thus, simplicity follows immediately from the previous corollary. These
algebras have real
rank zero because, as observed in \cite{Rordam}, the class of $\ca$-algebras of real rank zero is
closed under tensoring with finite dimensional $\ca$-algebras, forming direct sums, and forming
inductive limits. The Cuntz algebras $\O_{N^{n_k}}$ have real rank zero since they are purely
infinite. 
\bx

We finish this section by pointing out the connection between our results and the single variable
setting results of Bunce-Deddens. 

\begin{rem}
The focus of this paper is on the non-commutative multivariable
setting; however, we remark
that the proof of Theorem~\ref{periodk} goes through as presented for $N = 1$. Namely, the 
$\ca$-algebra generated by all  unilateral weighted shift operators of 
period $k$ is isomorphic
to $\M_k (\E_1)$, where $\E_1$ is the $\ca$-algebra generated by the 
unilateral shift, also realized as the algebra of 
Toeplitz operators with continuous symbol \cite{BD1,BD2}. The
proof presented here recaptures this result for $N =1$, although from a different
perspective. In particular, the Bunce-Deddens proof heavily relies on the associated function
theory which is omnipresent in the single variable case. Conceptually
speaking the proof here is more spatially oriented. 

While more effort is required to prove Theorem~\ref{periodk} for $N \geq 2$, 
the simplicity in Corollary 4.6 is more easily obtained as compared to the
single variable case. 
The basic point is that $\O_N$ is simple for $N\geq 2$, while for $N =1$ it is the         
$\ca$-algebra generated by the bilateral shift operator, the algebra of
continuous functions on the
unit circle, which is not simple. Nonetheless, the inductive limit
algebras $q (\fA (n_k))$
turn out to be simple for $N = 1$, and thus our results on these algebras 
for $N \geq 2$ can be
regarded as a non-commutative multivariable generalization of their 
result. 
\end{rem}

\section{Classification}\label{isos}

In this section, we establish an isomorphism theorem for the limit
algebras discussed in 
the previous two sections.
Let $\{ n_k \}_{k \geq 1}$ be an increasing sequence of positive integers with $n_k $ dividing
$n_{k+1}$ for $k \geq 1$.
Then for each prime $p$, there is a unique $\alpha_p$ in $\bbN \cup \{ \infty \}$ which is
the supremum of the exponents of $p$ which divide $n_k$ as $k \rightarrow \infty$. The 
{\it supernatural number} determined by the sequence $\{ n_k \}_{k \geq
1}$ is the formal
product $\delta (n_k) = \prod_{ p \, {\rm prime}} p^{\alpha_p}$. Given two such
sequences $\{ n_k \}_{k \geq 1}$ and $\{ m_j \}_{j \geq 1}$, it follows that $\delta (n_k) =
\delta (m_j)$ precisely when: for all $k \geq 1$, there is a $j\geq 1$
with $n_k | m_j$; and 
for all $j \geq 1$, there is a $k\geq 1$ with $m_j | n_k$. 

Supernatural numbers have been used to
classify UHF algebras \cite{Glimm1}, and   Bunce-Deddens algebras 
\cite{BD1,BD2}.
They also distinguish between the inductive limit algebras of the current 
paper, as do the 
associated $K_0$ groups. 

\begin{thm}\label{isothm}
Let $N \geq 2$ be a positive integer. Let $\{ n_k \}_{k \geq 1}$ and $\{ m_j \}_{j \geq 1}$ be
increasing sequences of positive integers for which $n_k | n_{k +1}$ and $m_j | m_{j +1}$ for
$j,k \geq 1$. Then the following are equivalent:
\begin{itemize}
\item[$(i)$] The supernatural numbers $\delta (n_k)$ and $\delta (m_j)$ are the same.
\item[$(ii)$] The algebras $\fA (n_k)$ and $\fA (m_j)$ are equal.
\item[$(iii)$] The algebras $q ( \fA (n_k) )$ and $q ( \fA (m_j))$ are equal.
\item[$(iv)$] If $\fB (n_k)$ is an inductive limit of Cuntz algebras determined by a sequence 
$B_{n_1} \rightarrow B_{n_2} \rightarrow \, \ldots$ such that $B_{n_k} \cong
\M_{d_{N,n_k}} ( \O_{N^{n_k}})$, and $\fB(m_j)$ is similarly defined, then $\fB(n_k)$ and
$\fB(m_j)$ are $*$-isomorphic. 
\item[$(v)$] The groups $K_0 ( \fB(n_k) )$ and $K_0 ( \fB(m_j) )$ are isomorphic.  
\end{itemize}
\end{thm}

\Prf
To see $(i) \Rightarrow (ii)$, observe the division property associated with $(i)$ shows that
each $n_k$ divides some $m_j$. Hence by Proposition~\ref{contain},
\[
\mathrm{C}^*_N(\mathrm{per}\,\,n_k) \subseteq \mathrm{C}^*_N(\mathrm{per}\,\,m_j)
\subseteq \fA(m_j)
\]
for $k \geq 1$. Whence, $\fA(n_k) \subseteq \fA(m_j)$. The converse inclusion follows by
symmetry.  The implications $(ii) \Rightarrow (iii)$ and $(iv) \Rightarrow (v)$ are obvious. 
Since $\fB (n_k) \cong q (\fA(n_k))$ and $\fB (m_j) \cong q (\fA(m_j))$ 
by  Corollary~\ref{qperk}, we have $(iii) \Rightarrow (iv)$. 

It remains to establish the implication $(v) \Rightarrow (i)$. Recall 
the $K_0$ group of $\O_N$ is the finite  abelian
group $K_0 (\O_N) = \bbZ / (N-1)\bbZ$ of order $N-1$. Hence, 
\[
K_0  (\M_{d_{N,n_k}}(\O_{N^{n_k}}) ) =  K_0  ( \O_{N^{n_k}}) = \bbZ / (N^{n_k} -1)\bbZ.  
\]
Since $\fB(n_k)$ is the inductive limit of the algebras $B_{n_k}$, we have
\[
K_0  ( \fB(n_k) ) = \lim_{\longrightarrow} K_0  ( B_{n_k} ).
\]
There are similar facts for $K_0  ( \fB(m_j) ) = \lim_{\longrightarrow} K_0  ( B_{m_j} ).$ 

For $k \geq 1$, let $g_k \in K_0  ( \fB(n_k) )$ be an element of order $N^{n_k} - 1$. Let 
$\Gamma : K_0  ( \fB(n_k) ) \rightarrow K_0  ( \fB(m_j) )$ be a group isomorphism. Then 
$\Gamma (g_k) \in K_0  ( \fB(m_j) )$, and it follows that the order of $\Gamma (g_k)$ must
divide $N^{m_j} - 1$ for some $j\geq 1$. (Every element of $K_0  (
\fB(m_j))$ has this 
property.)  

However, for positive integers $N \geq 2$ and $k\geq 1$, recall that
\[
d_{N,k} = 1 + N + \ldots + N^{k-1} = \frac{N^k -1}{N-1}. 
\]
Thus we have just observed that $d_{N,n_k}$ divides $d_{N,m_j}$. But this
implies that $n_k$ divides $m_j$. 
Indeed, suppose $d_{N,m_j} = c\, d_{N,n_k}$ for some positive integer
$c$. 
Consider the base $N^{n_k}$ expansion of $c$ given by 
$
c = c_0 + c_1 N^{n_k} + \ldots + c_l (N^{n_k})^l,
$
where $0 \leq c_i < N^{n_k}$ for $0 \leq i \leq l$. By comparing
coefficients in
$d_{N,m_j} = c \, d_{N,n_k}$, we get each $c_i = 1$ and $m_j -1= ln_k + 
n_k-1 $. Whence, $m_j =
(l+1)n_k$, and $n_k$ divides $m_j$.

By symmetry, every $m_j$ divides some $n_k$. Hence by the remarks
preceding the theorem, 
this shows that the supernatural numbers $\delta (n_k)$ and 
$\delta (m_j)$ are identical. 
\bx

\begin{rem}
After preparing this article, the author became aware of  related work
of Evans \cite{Evans} on Cuntz-Krieger algebras. The most notable overlap
between our papers is that the equivalence of conditions $(i)$ and $(v)$
in the previous theorem follows from the Cuntz algebra case of Theorem 4.3
from \cite{Evans}. We also point out that a related notion of periodicity
is used to define certain inductive limits of Cuntz-Krieger algebras in 
\cite{Evans}. In
the Cuntz algebra case, we see it is a more restrictive
version; requiring scalars $\{ c_k : k \geq 0\}$ such that $\lamiw =c_k$
for $1\leq i \leq N$ and all $|w| = k$. Thus the periodicity introduced
here is new, as is our main result Theorem 4.1. 
\end{rem}

We finish by pointing out that, not surprisingly, as for the Bunce-Deddens 
algebras, the 
algebras here are {\it not} almost finite dimensional. We need the
following easy generalization of a
theorem of Halmos. 

\begin{lem}
For $1\leq i \leq N$, the operator $L_i$ is not quasitriangular. 
\end{lem}

\Prf
In \cite{Halmos}, Halmos proves that the unilateral shift is not quasitriangular. However, as he
points out before proving this result, the proof really only depends on the operator of concern
being an isometry, with adjoint having non-trivial kernel. The $L_i$ 
clearly have this property
since $L_i^*$ annihilates the vacuum vector. 
\bx

We can follow the lines of the Bunce-Deddens proof to show these algebras are not AF.

\begin{thm}\label{af}
The algebras $ \fA (n_k) $ and $q( \fA (n_k))$ are not approximately finite dimensional. 
\end{thm}

\Prf
Suppose there are finite dimensional $\ca$-algebras $\{\fB_l\}_{l\geq 1}$ for which $\fB_l
\subseteq \fB_{l +1}$ and $q ( \fA(n_k) ) = \overline{\bigcup_{l\geq 1} \fB_l }.$ Then given
$\varepsilon > 0$, there is a $B_i \in \bigcup_{l\geq 1} \fB_l$ with $|| q(L_i) - B_i || <
\varepsilon$. Choose $A_i \in \fA (n_k)$ such that $q(A_i ) = B_i$. 
Then $|| q(L_i - A_i) || < \varepsilon$, so there is a compact operator $C_i$ with 
$|| L_i - A_i -C_i || < \varepsilon$. But $B_i$ belongs to a finite dimensional $\ca$-algebra,
hence there is a non-trivial polynomial $p$ with $p(B_i ) = p(q(A_i)) = 0$. Thus $A_i$ is
polynomially compact, and as such, it is quasitriangular. (This was proved
initially in
\cite{DouPea}.)   In particular, $A_i + C_i$ is quasitriangular, so that
$L_i$ belongs to the norm
closure of the quasitriangular operators, and is itself quasitriangular \cite{Halmos}. This
contradicts the previous lemma. Thus $q ( \fA(n_k) )$ is not approximately finite dimensional.
The proof that $\fA(n_k)$ is not AF is easier. 
\bx

{\noindent}{\bf Acknowledgements.}
The author is  grateful to Raul Curto, Ken Davidson, Palle Jorgensen and Paul Muhly
for enlightening conversations. Thanks also to the Department of Mathematics at the University
of Iowa for kind hospitality during the preparation of this article.




\end{document}